\documentclass{amsart}
\newtheorem{lemma}{Lemma}

\newcommand{\qu}{\mathbf{Q}}
\newcommand{\zet}{\mathbf{Z}}

\newtheorem{theorem}{Theorem}
\newtheorem{corollary}{Corollary}
\begin{document}
\title{On the distribution of Galois groups}
\author{Rainer Dietmann}
\address{Department of Mathematics, Royal Holloway, University of London\\
TW20 0EX Egham, United Kingdom}
\email{Rainer.Dietmann@rhul.ac.uk}
\maketitle
\begin{abstract}
Let $G$ be a subgroup of the symmetric group $S_n$, and let
$\delta_G=|S_n/G|^{-1}$ where $|S_n/G|$ is the index of $G$ in $S_n$.
Then there are at most
$O_{n, \epsilon}(H^{n-1+\delta_G+\epsilon})$
monic integer polynomials of degree $n$ having
Galois group $G$ and height not exceeding $H$,
so there are only `few' polynomials having `small' Galois group.
\end{abstract}
\section{Introduction}
On probabilistic grounds, one should expect that `almost all' polynomials
have the full symmetric group as Galois group acting on the roots.
More precisely, let
\begin{eqnarray*}
  E_n(H) & = & \#\{ (a_1, \ldots, a_n) \in \mathbf{Z}^n : |a_i| \le H
  \; (1 \le i \le n) \mbox{ and }\\
  & & \mbox{the splitting field of } X^n+a_1X^{n-1}+\ldots+a_n \\
  & & \mbox{ has not Galois group $S_n$}\}.
\end{eqnarray*}
Then van der Waerden (\cite{W}) showed that
\[
  E_n(H) \ll_n H^{n-6/((n-2)\log \log H)}.
\]
Using the large sieve, Gallagher (\cite{G}) improved this to
\begin{equation}
\label{der_hexer}
  E_n(H) \ll_n H^{n-1/2} \log H.
\end{equation}
Since there are asymptotically only
$c_n H^{n-1} + o(H^{n-1})$ polynomials
as above
which are reducible over $\mathbf{Q}$ (see \cite{C}),
one could conjecture that
$E_n(H) \ll_{n} H^{n-1}$.
Apart from an epsilon in the exponent,
this conjecture has been confirmed by Lefton (\cite{L}) for $n=3$
and by the author for $n=4$ (\cite{RD}). The following result gives
further support to this conjecture, showing that `small' Galois groups
are rare:
\begin{theorem}
\label{boeing_747}
Let $\epsilon>0$, and let
$n$ be a positive integer. Further, let $G$
be a subgroup of the symmetric group $S_n$, and let
\begin{eqnarray*}
  N_n(H; G) & = &
 \#\{ (a_1, \ldots, a_n) \in \mathbf{Z}^n : |a_i| \le H
  \; (1 \le i \le n) \mbox{ and } \nonumber \\
  & & X^n+a_1X^{n-1}+\ldots+a_n \mbox{ has Galois group $G$}\}.
\end{eqnarray*}
Then
\[
  N_n(H; G) \ll_{n,\epsilon} H^{n-1+\delta_G+\epsilon},
\]
where
\[
  \delta_G = \frac{1}{|S_n/G|}
\]
and $|S_n/G|$ is the index of $G$ in $S_n$.
\end{theorem}
Note that in particular we immediately obtain the bound $E_n(H)
\ll_{n, \epsilon} H^{n-1/2+\epsilon}$ from Theorem
\ref{boeing_747}, which is only slightly worse
than Gallagher's bound (\ref{der_hexer}). However, if we exclude
the worst case $|S_n/G|=2$, we can do much better. This is a
situation occurring sometimes in applications, for instance
Zarhin \cite{Z} proved that the Galois action on the Tate modules
of a Jacobian variety of a curve $C/\qu$ of the type $y^2=f(x)$
is `big' provided deg$(f) \ge 5$ and the Galois group of $f$
contains the alternating group $A_n$ on $n$ letters.
So let
\begin{eqnarray*}
  E_n(H)' & = & \#\{ (a_1, \ldots, a_n) \in \mathbf{Z}^n : |a_i| \le H
  \; (1 \le i \le n) \mbox{ and }\\
  & & X^n+a_1X^{n-1}+\ldots+a_n \mbox{ has not Galois group $S_n$ or $A_n$}\}.
\end{eqnarray*}
Clearly the upper bound in (\ref{der_hexer}) also holds true for the
quantity $E_n(H)'$, but we can improve this considerably by
using the fact that transitive proper
subgroups of $S_n$ different from $A_n$ have large index in terms of $n$:
\begin{corollary}
\label{como}
For $n \ge 9$ we have
\[
 E_n(H)' \ll_{n, \epsilon} H^{n-1+e(n)+\epsilon}
\]
where
\[
  e(n) = \frac{2}{{n \choose {\lfloor n/2 \rfloor}}}
\]
and $\lfloor \frac{n}{2} \rfloor$ is the largest integer not exceeding
$\frac{n}{2}$.
\end{corollary}
Note that $e(n)$ tends rapidly to zero as $n$ approaches infinity.
Since by Chela's result \cite{C} there are only $O(H^{n-1})$ reducible
polynomials $f$ in question, and since irreducible $f$ have a splitting
field with transitive Galois group, Corollary \ref{como} follows
immediately from Theorem \ref{boeing_747} and Lemma \ref{boeing_767} in
\S \ref{y}.\\
For even degree trinomials, we can also handle $A_n$
and obtain a bound which for large $n$ is quite close to the expected one:
\begin{theorem}
\label{montreaux}
Let $n$ and $r$ be coprime positive integers where $n>r$, $n \ge 10$
and $n$ is even. Then
\begin{eqnarray*}
  \#\{ (a_r, a_n) \in \mathbf{Z}^2 : |a_r|, |a_n| \le H
  \mbox{ and }\\
  X^n+a_rX^{n-r}+a_n \mbox{ has not Galois group $S_n$}\}
  \ll_{n, \epsilon} H^{1+\frac{1}{n} +\epsilon}.
\end{eqnarray*}
\end{theorem}
Note that the condition $gcd(r,n)=1$ here is obviously necessary
(otherwise the Galois group would always be strictly smaller
than $S_n$), whereas
forcing $n$ to be even is a restriction resulting from the method of proof.\\ \\
Unlike previous approaches to this kind of problem using sieve methods
incorporating local information modulo $p$ for many primes $p$, we work
globally, using Galois resolvents and this way reduce the original problem
to one about counting integer points on certain varieties. This is the
reason why we are able to go beyond a square root saving, which is typically
the best one can obtain by the large sieve. However, one should note that
one of our key ingredients (see \cite{HB}) also makes use of suitable reductions
modulo $p$.\\
Since Gallagher's paper \cite{G} has had a few applications and extensions
since its publication, for example to polynomials with some of the
coefficients $a_i$ fixed (\cite{Coh}), reciprocal polynomials
(\cite{DDS}), characteristic polynomials of unimodular matrices
(\cite{R}) or $L$-functions of algebraic curves (\cite{Ch}, \cite{K}),
one might wonder if our methods are also applicable in those cases.
However, it is not obvious how to extend our methods to
those problems. One of the main obstacles seems to be a suitable analogue
of Lemma \ref{kitchener}. \\ \\
\textit{Acknowledgment.} The author would like to thank Dr. T.D. Browning
for bringing the reference \cite{BH} to his attention.
\section{Preparations}
\label{y}
\begin{lemma}
\label{kapstadt}
Let $f(X)=a_0X^n+a_1X^{n-1}+\ldots+a_n \in \mathbf{C}[X]$. Then all roots
$z \in \mathbf{C}$ of the equation $f(z)=0$ satisfy the inequality
\[
|z| \le \frac{1}{\sqrt[n]{2}-1} \cdot \max_{1 \le k \le n}
\sqrt[k]{\left| \frac{a_{k}}{a_0 {n \choose k}} \right|}.
\]
\end{lemma}
\proof
This is Theorem 3 in \S27 of \cite{M}.
\begin{lemma}
\label{kitchener}
Let $n$ and $r$ be coprime positive integers where $n>r$,
and let $a_1, \ldots, a_{r-1},$
$a_{r+1}, \ldots, a_{n-1}$ be fixed integers. Then the polynomial
\[
X^n + a_1 X^{n-1} + \ldots + a_{n-1} X + t 
\]
has for all but at most $n^2+n$ integers $a_r$ the full symmetric group $S_n$
as Galois group over the rational function field $\qu(t)$.
\end{lemma}
\proof
This follows from the proof of Satz 1 in \cite{H1}; see also
the introduction of \cite{H2}.

\begin{lemma}
\label{boeing_767}
Let $n$ be a positive integer with $n \ge 9$, and let $G$
be a transitive subgroup of $S_n$ with $G \ne S_n$ and $G \ne A_n$. Then
\[
|S_n/G| \ge \frac{1}{2} {{n} \choose {\lfloor \frac{n}{2}} \rfloor}.
\]
\end{lemma}
\proof
This follows immediately from Theorem 5.2B in \cite{DM} (using
$r=\lfloor \frac{n}{2} \rfloor$).
\begin{lemma}
\label{toronto}
Let $n$ be a positive integer with $n \ge 9$, and let $G$ be a subgroup
of $S_n$ different from $S_n$ and $A_n$. Then either $G$ has index $n$
or index $2n$
in $S_n$, and leaves one element in $\{1, \ldots, n\}$ fixed,
or $G$ has index at least $\frac{n(n-1)}{2}$ in $S_n$.
\end{lemma}
\proof
This follows immediately from Theorem 5.2B in \cite{DM} (using $r=2$).
\begin{lemma}
\label{resolv}
Let $n$ be a positive integer, $G$ a subgroup of $S_n$,
\[
f(X) = X^n + a_1 X^{n-1} + \ldots + a_n \in \zet[X]
\]
with complex roots $\alpha_1, \ldots, \alpha_n$ and
\[
\Phi(z; a_1, \ldots, a_n) = \prod_{\sigma \in S_n/G} \left(
z-\sum_{\tau \in G} \alpha_{\sigma(\tau(1))}
\alpha_{\sigma(\tau(2))}^2 \cdots \alpha_{\sigma(\tau(n))}^n \right).
\]
Then this \textit{Galois resolvent} $\Phi(z; a_1, \ldots, a_n)$
is a polynomial in $z, a_1, \ldots, a_n$ with integer coefficients.
Moreover,
if the splitting field of $f$ over the rationals has Galois group $G$,
then $\Phi(z; a_1, \ldots, a_n)$ has a rational
(and thus integral) root $z$.
\end{lemma}
\proof
This is a well known result, see for example Lemma 3.2 in \cite{KM}.
\begin{lemma}
\label{athen}
Let $F \in \zet[X_1, X_2]$ be of degree $d$
and irreducible over $\qu$. Further,
let $P_1, P_2$ be real numbers such that $P_1, P_2 \ge 1$, and let
\[
  N(F; P_1, P_2) = \#\{\mathbf{x} \in \zet^2:
  F(\mathbf{x})=0
  \mbox{ and } |x_i| \le P_i \; (1 \le i \le 2)\}.
\]
Moreover, let
\[
  T = \max \left\{ \prod_{i=1}^2 P_i^{e_i} \right\}
\]
with the maximum taken over all integer $2$-tuples $(e_1, e_2)$ for
which the corresponding monomial $X_1^{e_1} X_2^{e_2}$ occurs
in $F(X_1, X_2)$ with nonzero coefficient. Then
\begin{equation}
\label{pigs}
  N(F; P_1, P_2) \ll_{d, \epsilon}
  \max\{P_1, P_2\}^{\epsilon}
  \exp \left( \frac{\log P_1 \log P_2}{\log T} \right).
\end{equation}
\end{lemma}
\proof
After homogenizing $F$, this becomes
the special case $P_1=1$ of Theorem 1 in \cite{BH}; see also
formula (3) there, and see \cite{HB2}, Theorem 15 for an earlier reference for
this result.
Note that in \cite{BH}, it is assumed that $F$ is
absolutely irreducible rather than irreducible over the rationals, but
as remarked in \cite{HB}, top of page 556, if $F$ is irreducible over
the rationals, but not absolutely irreducible, then by B\'{e}zout's
Theorem $N(F; P_1, P_2) \ll_d 1$, so we may clearly assume that
$F$ is absolutely irreducible.\\ \\
\textit{Remark.} 
Bombieri and Pila's result \cite{BP} is essentially the special
case $P_1=P_2$ of Lemma \ref{athen}.
Since in our later application of this lemma
our box will be rather
lopsided, we need this more powerful generalization.

\section{Proof of Theorem \ref{boeing_747}}
To prove Theorem \ref{boeing_747}, we have to show that if $G$ is a
subgroup of $S_n$ of index $m$, then 
\begin{eqnarray}
\label{englisch}
\#\{ (a_1, \ldots, a_n) \in \mathbf{Z}^n : |a_i| \le H
\; (1 \le i \le n) \nonumber \\ \mbox{ and }
X^n+a_1X^{n-1}+\ldots+a_n \mbox{ has Galois group $G$}\}
\ll_{n,\epsilon} H^{n-1+\frac{1}{m}+\epsilon}.
\end{eqnarray}
So fix such $G$ and let $\Phi(z; a_1, \ldots, a_n)$
be the corresponding Galois resolvent from Lemma \ref{resolv}.
By Lemma \ref{resolv}, the resolvent
$\Phi(z; a_1, \ldots, a_n)$ is of the form
\begin{equation}
\label{ulm}
\Phi(z; a_1, \ldots, a_n) = z^m + b_1(a_1, \ldots, a_n) z^{m-1} +
\ldots + b_m(a_1, \ldots, a_n),
\end{equation}
where the $b_i$ are suitable polynomials in $a_1, \ldots, a_n$ with
integer coefficients. Now if
\begin{equation}
\label{portugal}
|a_i| \le H \;\;\; (1 \le i \le n)
\end{equation}
and
\[
f(X)=X^n+a_1 X^{n-1} + \ldots + a_n
\]
has Galois group $G$, then by Lemma \ref{resolv} the resolvent
$\Phi(z)=\Phi(z; a_1,
\ldots, a_n)$ has an integer root $z$. We seek to obtain an upper bound
on $|z|$. From Lemma \ref{kapstadt} it is clear that there exists a positive
constant $\alpha$ depending at most on $n$ such that $\alpha \ge 1$ and
$|z| \ll H^{\alpha}$. Here and in the following  we adopt the convention
that all implied constants depend at most on $n$ and $\epsilon$.
Thus if (\ref{portugal}) holds and $f$
has Galois group $G$, then $\Phi(z)$ has an integer root $z$ with
$|z| \ll H^{\alpha}$.\\
Now choose $a_1, \ldots, a_{n-2} \in \mathbf{Z}$ with
$|a_i| \le H \; (1 \le i \le n-2)$. We want to bound the number of integers
$a_{n-1}, a_n$ with
$|a_{n-1}|, |a_n| \le H$ such that $f$ has Galois group $G$. With respect
to (\ref{englisch}), it is sufficient to show that there are at most
$O(H^{1+\frac{1}{m}+\epsilon})$
such $a_{n-1}, a_n$. By Lemma \ref{kitchener},
for our fixed $a_1, \ldots, a_{n-2}$ the polynomial
\[
X^n + a_1 X^{n-1} + \ldots + a_{n-1} X + t
\]
has for all but at most $n^2+n$
values of $a_{n-1}$ the full symmetric group $S_n$ as Galois group over
the rational function field $\mathbf{Q}(t)$. With respect to
(\ref{englisch}) it thus suffices to fix any such $a_{n-1}$ with
$|a_{n-1}| \le H$ for which
\[
X^n + a_1 X^{n-1} + \ldots + a_{n-1} X + t
\]
has Galois group $S_n$ over $\mathbf{Q}(t)$ and then show that for those
fixed $a_1, \ldots, a_{n-1}$ there are at most
$O(H^{\frac{1}{m}+\epsilon})$
possible $a_n \in \mathbf{Z}$ with $|a_n| \le H$ for which
\[
f(X)= X^n + a_1 X^{n-1} + \ldots + a_n
\]
has Galois group $G$.\\
So let us consider
the Galois resolvent $\Phi(z, a_n) = \Phi(z; a_1, \ldots, a_n)$ as a
polynomial in $z$ and $a_n$. Since
\[
X^n + a_1 X^{n-1} + \ldots +
a_{n-1} X + t
\]
has Galois group $S_n$ over $\mathbf{Q}(t)$, the resolvent
$\Phi(z, a_n)$ must be irreducible over the rationals. Moreover,
by (\ref{ulm}) the resolvent
$\Phi(z, a_n)$ obviously has degree at least $m$ in $z$,
since an additive term $z^m$ occurs. Further, if $f$ has Galois group
$G$ for a particular choice of $a_n$ with $|a_n| \le H$,
then $\Phi(z)=\Phi(z, a_n)$ must
have an integer root $z$ with $|z| \ll H^{\alpha}$. So it is
sufficient to bound above the number of integer zeros of $\Phi(z, a_n)$
with $|z| \ll H^{\alpha}$ and $|a_n| \le H$.
To this end
we apply Lemma \ref{athen} with $P_1 \asymp H^\alpha$ and $P_2 = H$.
Since $\Phi(z, a_n)$ has degree at least $m$ in $z$,
we have $T \gg H^{m \alpha}$.
Hence by Lemma \ref{athen} the quantity
\[
  \#\{(z, a_n) \in \mathbf{Z}^2 : |z| \ll H^{\alpha}, |a_n| \le H
  \mbox{ and } \Phi(z,a_n)=0\} \\
\]
under scrutiny can be bounded by
\[
  H^{\epsilon}
  \exp \left( \frac{\log H}{m} + O(1) \right) \ll H^{1/m+\epsilon}.
\]
This completes the proof of the theorem.

\section{Proof of Theorem \ref{montreaux}}
Let us first deal with irreducible $f(X)=X^n+a_rX^{n-r}+a_n$. Then
the splitting field of $f$ has a transitive Galois group $G$.
Let us first assume that $G \ne A_n$, $G \ne S_n$. Then we can just
copy the proof of Theorem \ref{boeing_747}, checking that still $f$
generically has Galois group $S_n$ by Lemma \ref{kitchener}.
Lemma \ref{boeing_767} then produces a saving good
enough, because $n \ge 10$. If $G=A_n$, then the reasoning in section $3$
of \cite{L} gives what we want, since $n$ is even: Using an
explicit discriminant formula being quadratic in $a_r$ and
reducing the problem to counting points on conics one can even
save $-1+\epsilon$ in the exponent. Reducible $f$
can be handled the following way: Suppose first that such $f$ has a rational
and thus integer zero $x$.
The number of such $f$ can be bounded
as follows: By Lemma \ref{kapstadt} we have $|x| \ll H$, and
by Lemma \ref{kitchener},
for all but at most $n^2+n$ choices of $a_r$ the polynomial
$X^n+a_rX^{n-r}+a_n$ as a polynomial in $X$ and $a_n$ is irreducible.
Now Lemma \ref{athen} gives at most
$O(H^{\frac{1}{n}+\epsilon})$ choices for $a_n$ and $x$ with
$|a_n| \le H, |x| \ll H$, and by taking the $a_r$ into account
we find at most $O(H^{1+\frac{1}{n}+\epsilon})$ such trinomials $f$
having an integer zero.
If $f$ is reducible, but has no integer zero, then
$f$ splits off in factors each of degree at least two over the rationals.
Let $G$ be the Galois group of the splitting field of $f$. Then $G$ must
be intransitive, hence different from $S_n$ and $A_n$, and $G$ cannot
leave one element fixed, since
otherwise $f$ would split off a linear factor, implying an integer zero
for $f$. By Lemma \ref{toronto} we conclude that $G$ has index at least
$\frac{n(n-1)}{2}$ in $S_n$ and again we can copy the proof of Theorem
\ref{boeing_747} for this $G$ and trinomials $f$, obtaining the bound
$O(H^{1+\frac{1}{n(n-1)/2}+\epsilon})= O(H^{1+\frac{1}{n}+\epsilon})$
since $n \ge 10$.
This completes the proof of Theorem \ref{montreaux}.


\begin{thebibliography}{10}
\bibitem{BP} \textsc{Bombieri, E. \& Pila, J.}
{The number of integral points on arcs and ovals,
\textit{Duke Math. J.} \textbf{59} (1989), 337--357.}
\bibitem{BH} \textsc{Browning, T.D. \& Heath-Brown, D.R.}
{Plane curves in boxes and equal sums of two powers,
\textit{Math. Z.} \textbf{251} (2005), 233--247.}
\bibitem{Ch} \textsc{Chavdarov, N.}
{The generic irreducibility of the numerator of the zeta function in a
family of curves with large monodromy, 
\textit{Duke Math. J.} \textbf{87} (1997), 151--180.}
\bibitem{C} \textsc{Chela, R.} {Reducible polynomials,
\textit{J. London Math. Soc.} \textbf{38} (1963), 183--188.}
\bibitem{Coh} \textsc{Cohen, S.D.} {The distribution of the Galois groups of
integral polynomials, \textit{Illinois J. Math.} \textbf{23} (1979),
135--152.}
\bibitem{DDS} \textsc{Davis, S., Duke, W. \& Sun, X.}
{Probabilistic Galois theory of reciprocal polynomials,
\textit{Exposition. Math.} \textbf{16} (1998), 263--270.}
\bibitem{RD} \textsc{Dietmann, R.} {Probabilistic Galois theory for
quartic polynomials, \textit{Glasgow Math. J.} \textbf{48} (2006),
no. 3, 553--556.}
\bibitem{DM} \textsc{Dixon, J.D. \&  Mortimer, B.}
{Permutation groups, Graduate Texts in Mathematics,
Springer-Verlag, New York (1996).}
\bibitem{G} \textsc{Gallagher, P.X.} {The large sieve and probabilistic
Galois theory, \textit{Proceedings of Symposia in Pure Mathematics} XXIII
(1973, A.M.S.)}
\bibitem{HB} \textsc{Heath-Brown, D.R.} {The density of rational points
on curves and surfaces, \textit{Ann. of Math.} \textbf{155} (2002),
553--598.}
\bibitem{HB2} \textsc{Heath-Brown, D.R.} {Counting rational points on
algebraic varieties, Springer Lecture Notes \textbf{1891} (2006), 51--95.}
\bibitem{H1} \textsc{Hering, H.} {Seltenheit der Gleichungen mit Affekt
bei linearem Parameter, \textit{Math. Ann.} \textbf{186} (1970), 263--270.}
\bibitem{H2} \textsc{Hering, H.} {\"Uber Koeffizientenbeschr\"ankungen
affektloser Gleichungen, \textit{Math. Ann.} \textbf{195} (1972), 121--136.}
\bibitem{KM} \textsc{Kl\"uners, J. \& Malle, G.} {Explicit Galois
realization of transitive groups of degree up to $15$,
\textit{J. Symbol. Comp.} \textbf{30} (2000), 675--716.}
\bibitem{K} \textsc{Kowalski, E.}
{The large sieve, monodromy and zeta functions of curves,
\textit{J. Reine Angew. Math.} \textbf{601} (2006), 29--69.}
\bibitem{L} \textsc{Lefton, P.} {On the Galois groups of cubics
and trinomials, \textit{Acta Arith.} \textbf{XXXV} (1979), 239--246.}
\bibitem{M} \textsc{Marden, M.} {Geometry of polynomials,
second edition, Mathematical Surveys, No. 3, \textit{American
Mathematical Society} (1966).}
\bibitem{R} \textsc{Rivin, I.}
{Walks on groups, counting reducible matrices, polynomials, and surface and
free group automorphisms, \textit{Duke Math. J.} \textbf{142} (2008), 353--379.}
\bibitem{W} \textsc{van der Waerden, B.L.} {Die Seltenheit der reduziblen
Gleichungen und die Gleichungen mit Affekt, \textit{Monatsh. Math.}
\textbf{43} (1936), 137--147.}
\bibitem{Z} \textsc{Zarhin, Y.G.}
{Very simple $2$-adic representations and hyperelliptic Jacobians,
\textit{Mosc. Math. J.} \textbf{2} (2002), 403--451.}
\end{thebibliography}
\end{document}